\documentstyle [11pt,graphicx,amstex,epsfig,amscd,amssymb]{article}
\textwidth
6.2in

\def\om{\omega}

\def\phi{\varphi}

\def\g{\mathfrak g}
\def\l{\lambda}

\def\ti{\times}

\def\sof{{\mathfrak so}}

\def\be{\begin{equation}}
\def\ee{\end{equation}}
\def\bear{\begin{eqnarray}}
\def\eear{\end{eqnarray}}
\def\best{\begin{eqnarray*}}
\def\eest{\end{eqnarray*}}
\def\pf{{ {\emph{Proof}}}: }
\def\ni{\noindent}

\newtheorem{theorem}{Theorem}[section]

\newtheorem{proposition}[theorem]{Proposition}

\newtheorem{defn}[theorem]{Definition}

\newtheorem{remark}[theorem]{Remark}

\newtheorem{example}[theorem]{Example}

\renewcommand{\Re}{\operatorname{Re}}

\renewcommand{\Im}{\operatorname{Im}}

\def\r#1{\right#1}
\def\l#1{\left#1}

\def\Z{{ \Bbb Z}}
\def\R{{ \Bbb R}}
\def\C{{ \Bbb C}}

\def\ts3{{T^*S^3}}

\title{\bf  Special Lagrangians of cohomogeneity one in the resolved conifold}

\author{Marianty Ionel and Maung Min-Oo}

\date{}

\begin{document}


\maketitle

\vskip.1in

\noindent {\bf Abstract.} In this paper, which is a natural continuation of \cite{im1}, we describe some special Lagrangians of cohomogeneity one in the resolved conifold. 
Our main result gives a foliation of the resolved conifold by $T^2$-invariant special Lagrangians, where the generic leaf is topologically $T^2 \times \R$. We also obtain a family of  $SO(3)$-invariant special Lagrangians. These special Lagrangian families in both the deformed and the resolved conifold approach asymptotically the same special Lagrangian cones in the conifold.
\\

\\



\section{Introduction}

The notion of a special Lagrangian submanifold was introduced in 1984 by Harvey and Lawson in their seminal paper \cite{hl} on calibrated geometries.
Special Lagrangian geometry in Calabi-Yau manifolds has become an
important subject due to a phenomenon in physics known
as mirror symmetry. In 1996, Strominger, Yau and Zaslow
\cite{syz} conjectured that a compact Calabi-Yau $3$-fold and its
mirror should be foliated by special Lagrangian $3$-tori with
possibly singular fibres and the fibrations are dual to each
other. This conjecture proposes a way to construct the mirror of a
compact Calabi-Yau manifold, by an appropriate compactification of
the dual of the special Lagrangian fibration. Since compact
Calabi-Yau metrics cannot be written explicitly and also
because singularities play a fundamental role, physicists have always been interested in finding concrete 
Calabi-Yau metrics on non-compact manifolds. One of the earliest examples of a pair of mirror Calabi-Yau metrics was found by Candelas and de la Ossa \cite{co} in 1990. The two manifolds arise from perturbing a singular cone on $S^2 \times S^3$ and are respectively known as the deformed and the resolved conifold. In fact, there is a $1$-parameter family of Calabi-Yau metrics that passes through the singular metric and transforms the deformed conifold into the resolved conifold. The deformed conifold is a (trivial) $\R^3$-bundle over $S^3$ and the resolved conifold is a $\C^2$-bundle over $S^2$. In passing through the singularity the special Lagrangian $S^3$ in the deformed conifold is pinched to a point and reappears in the resolved conifold as a holomorphic $\C P^1$. This is known as the {\it conifold transition} (see  \cite{co} and \cite {sty}).

In our recent paper \cite{im1} we obtained a description of
cohomogeneity one special Lagrangians
 in the deformed conifold. In particular, we found that the deformed conifold is foliated by special Lagrangians
  which are generically $T^2 \times \R$. In the present paper we investigate the analogous problem for
  the resolved conifold and exhibit a foliation of the resolved conifold by $T^2$-invariant special Lagrangians, where the generic leaf is topologically $T^2 \times \R$. 
We also obtain an $SO(3)$-invariant family of special Lagrangians.
It turns out that in both cases, the special Lagrangian families
found in the resolved conifold have the same asymptotics as the
special Lagrangians that we described in the
deformed conifold in \cite{im1}, i.e. both families asymptotically
approach the same special Lagrangians cones in the conifold.
\bigskip


\section{Special Lagrangian Geometry}

 In what follows, we will review some definitions and
set up notation.

\begin{defn}
A Calabi-Yau $n$-fold $(M,J,\omega,\Omega)$ is a K{\"a}hler
$n$-dimensional manifold $(M,J,\omega)$ with Ricci-flat K{\"a}hler
metric $g$ and a nonzero holomorphic section $\Omega$ which
trivializes the canonical bundle $K_M$.
\end{defn}

Since the metric $g$ is Ricci-flat, $\Omega$ is a parallel tensor
with respect to the Levi-Civita connection $\nabla^g$ \cite{dj1}.
By rescaling $\Omega$, we can take it to be the holomorphic
$(n,0)$-form that satisfies:
\begin{equation}
{\frac{\om^n}{n!}}={(-1)^{\frac{n(n-1)}{2}}{\biggl(\frac{i}{2}\biggl)}^n\Omega\wedge\bar
\Omega},
\end{equation}
where $\om$ is the K{\"a}hler form of $g$. The form $\Omega$ is
called the {\it holomorphic volume form} of the Calabi-Yau
manifold $M$.
Such a holomorphic volume form $\Omega$ exists if and only if the
canonical bundle $K_M$ is trivial, which implies that the first
Chern class $c_1(M)=0$ in $H^2(M,\Z)$. Also, the above definition
is equivalent to the manifold $(M,g)$ having holonomy $SU(n)$.

\begin{defn}
Let $M^n$ be an $n$-dimensional Calabi-Yau manifold and $L\subset
M$ a real $n$-dimensional submanifold of $M$. Then $L$ is a {\it
special Lagrangian submanifold} of $M$ if $ \om\mid_L\equiv 0$
(i.e. $L$ is Lagrangian) and $\Im\Omega\mid_L\equiv 0$.
Equivalently, $\Re(\Omega)\mid_L$ is the volume form on $L$ with
respect to the induced metric.
\end{defn}

One important property of the special Lagrangian submanifolds is
that they are absolutely area minimizing in their homology class,
so they are in particular {\it minimal} submanifolds (see
\cite{hl}).

 \noindent An example of a Calabi-Yau manifold
is $\C^n$, with coordinates $(z_1,\dots,z_n)$, endowed with the
flat metric $g$, the K{\"a}hler form $\omega_0$ and the
holomorphic volume form $\Omega_0=dz_1\wedge\dots\wedge dz_n $.
\noindent $\R^n$ is a trivial example of a special Lagrangian
submanifold in $\C^n$. The following standard examples of
cohomogeneity one special Lagrangians in $\C^n$ are due to
Harvey-Lawson.
\smallskip

{\noindent{\bf Example 1:}
\begin{equation}
\label{eqHL}
 L_c= \{ \lambda \, u  \,\,|\, \, u\in S^{n-1}\subset\R^n, \,\, \lambda \in \C,
\ \Im(\lambda^n)=c \, \},
\end{equation}
\ni are special Lagrangians submanifolds of $\C^n$, invariant
under the action of $SO(n)\subset SU(n)$, where $c\in \R$ is a
constant.

\ni 
{\it Remark:} The variety $L_0$ is an union of $n$ special
Lagrangian $n$-planes and when $c\not = 0$, each component of
$L_c$ is diffeomorphic to $\R\ti S^{n-1}$ and it is asymptotic to
$L_0$.

\medskip
\noindent {\bf Example 2:} The special Lagrangian submanifolds in
$\C^n$, invariant under the action of the maximal torus
$T^{n-1}=\{\mbox{ diag}(e^{i\theta_1}, \dots, e^{i\theta_n})|
\theta_1+ \dots \theta_n=0\}$ of $SU(n)$ are given by the
equations:

\begin{align*}
&|z_1|^2-|z_j|^2=c_j, \,\,j=2,3,\dots,n \mbox{ \quad    and}\\
& \Re (z_1 z_2 \dots z_n)=c_1,\, \mbox {  if $n$ even}\\
& \Im (z_1 z_2 \dots z_n)=c_1,\, \mbox {  if $n$ odd}
\end{align*}
where $c_1,c_2,\dots c_n$ are any real constants.
\smallskip

\ni {\it Remark:} These special Lagrangians are topologically
$T^{n-1} \times \R\,$ and they are asymptotic to
cones on two flat tori in $S^{2n-1}$ obtained by setting all
constants equal to $0$.


\section{The resolved conifold as a Calabi-Yau manifold}

In \cite{im1}, we studied the cohomogeneity one special Lagrangian
submanifolds in the deformed conifold $\ts3$. In this paper we
will answer the same question in the resolved conifold. Both manifolds
are Calabi-Yau and asymptotic to the cone on $S^2\times S^3$ with
a Ricci flat metric, called the {\it conifold}. Following the original paper of Candelas-de la Ossa \cite{co}, we will first give a description of the resolved conifold and its Calabi-Yau structure.

\medskip
\ni 
Let $Q_0$ be the singular quadric in $\C^{4}$ defined by the
equation:
\begin{equation}\label{conif}
\sum_{i=0}^{3}z_i^2=0
\end{equation}
This quadric, called the conifold, is a cone on $S^2\times S^3$. The deformed conifold is obtained by perturbing the conifold equation to $\sum_{i=0}^{3}z_i^2=\epsilon^2$, with $\epsilon$ a positive constant,
and it admits a Calabi-Yau structure.

\ni The other way to deal with the singularity of the conifold is
to resolve it. Making the linear change of variables:

\begin{equation}\label{coord}
\begin{pmatrix} X\\Y\\U\\V\end{pmatrix}= \frac{1}{\sqrt{2}}\begin{pmatrix} 1&-i&0&0\\1&i&0&0\\0&0&-i&1\\0&0&-i&-1
\end{pmatrix}\begin{pmatrix} z_0\\z_1\\z_2\\z_3\end{pmatrix}
\end{equation}

\ni
transforms the quadratic form $\, \sum_{i=0}^{3}z_i^2\, $ to $\, 2(XY-UV)\,$
and the conifold equation (\ref{conif}) becomes:

\begin{equation}\label{coneq}
XY-UV=0
\end{equation}

\ni We note that the matrix
\begin{equation}\label{matrix}
P=\frac{1}{\sqrt{2}}\begin{pmatrix}1&-i&0&0\\1&i&0&0\\0&0&-i&1\\0&0&-i&-1\end{pmatrix}
\end{equation}
is not in $SO(4,\C)$ (because it transforms the standard quadratic
form), but it is an element of $U(4)$, since $P^{-1}= P^* = {\bar
P}^t$. This shows that the Euclidean length is still preserved:

\begin{equation}\label{length}
r^2 = \sum_{i=0}^{3}|z_i|^2\ = |X|^2 +|Y|^2 +|U|^2 +|V|^2
\end{equation}

\ni The $SO(4,\C)$ action on the $z$ variables is now conjugated
to an action on the new variables via $g \mapsto {\tilde g} =
PgP^{-1}$,  i.e.,
\begin{equation}\label{action}
{\tilde g}\begin{pmatrix}X\\Y\\U\\V \end{pmatrix} =
P\,g\,P^*\begin{pmatrix}X\\Y\\U\\V \end{pmatrix}, \ g \in SO(4,\C)
\end{equation}
\ni

\ni
Let
\begin{equation}\label{W}
W=\l(\begin{matrix}X&U\\V&Y\end{matrix}\r)=
\frac{1}{\sqrt{2}}\l(\begin{matrix} z_0-iz_1 & z_3-iz_2 \\ -z_3-iz_2 & z_0+iz_1
\end{matrix}\r)
\end{equation}

\ni
To resolve the conifold, one has to replace the equation
$(\ref{coneq})$ by the pair of equations:

\begin{equation}\label{resconif}
\l(\begin{matrix}X&U\\V&Y\end{matrix}\r)\l(\begin{matrix}\lambda_1\\\lambda_2\end{matrix}\r)
=\l(\begin{matrix}0\\0\end{matrix}\r)
\end{equation}

\ni
where $[\lambda_1:\lambda_2]$ are homogeneous coordinates on $\C
P^1 = S^2$. In other words, the resolved conifold $M$ is the complex 3-fold defined by:

\begin{equation}\label{rescon}
\l \{ (X,Y,U,V,[\lambda_1:\lambda_2])\in \C^4 \times \C P^1 \,|\, XY - UV = 0,\, X\lambda_1 + U \lambda_2 = 0,
V\lambda_1 + Y \lambda_2 = 0 \r \}
\end{equation}

\ni Each point $(X,Y,U,V)\in Q_0$,  except the origin, determines
a unique point
 $$[\lambda_1:\lambda_2]=  [-U:X] = [-Y:V] \in \C P^1$$ 

\ni
and the singularity at the origin is replaced (blown-up) by a copy of $\C  P^1$. Outside the origin, the resolved conifold is topologically the same as the conifold.

\medskip
\ni
{\it Remark:} 
There is a dual way to resolve the singularity where we replace the equation (\ref{resconif}) by its transpose
$$\begin{pmatrix}\lambda_1 & \lambda_2\end{pmatrix} \begin{pmatrix}X&U\\V&Y\end{pmatrix} =\begin{pmatrix}0 &
0\end{pmatrix}
$$
These are known as small resolutions of the conifold. To obtain
the ``full'' resolution of the singularity, one replaces the
origin with the set of all complex lines through the origin that
lie on the cone and thus the singularity would be replaced by an
$S^2 \times S^2$ rather than a single $S^2$.

\medskip
\ni
We will also use inhomogeneous coordinates $\lambda_+
=\frac{\lambda_2}{\lambda_1}$ in the coordinate patch $H_+$ where
$\lambda_1\not=0$ and $\lambda_{-}=\frac{\lambda_1}{\lambda_2}$ in
the coordinate patch $H_{-}$, where $\lambda_2\not=0$. In $H_+$ we
can take $(U,Y, \lambda_+)$ and in $H_{-}$ we can take $(X,V,
\lambda_{-})$ as coordinates for the conifold. On the intersection
$H_+\cap H_-$ we have:
$$(X,V,\lambda_{-})=(-\lambda_+ U,-\lambda_+ Y,\lambda_+^{-1})$$

\ni
In \cite{co}, Candelas and de la Ossa showed that the resolved conifold is in fact the
${\mathcal O}(-1)\oplus {\mathcal O}(-1)$ bundle over $\C P^1\simeq S^2$ with fibre $\C^2$. 
The radius in the fibre is given by
$$r^2=\mbox{tr}(W^*W)=(1+|\lambda_+|^2)(|U|^2+|Y|^2)=(1+|\lambda_{-}|^2)(|X|^2+|V|^2)$$
Asymptotically, as $r\to \infty$, the resolved conifold approaches the singular conifold $Q_0$.

\ni 
The Ricci-flat metric on $M$ found in \cite{co} is:
\begin{equation} \label{metric}
g_{rc}=F_{a}'(r^2)\,\mbox{tr}(dW^*dW) + F_{a}''(r^2)|\mbox{tr}(W^*dW)|^2+4a^2 g_{S^2}
\end{equation}
where $F_{a}$ is a function of $r^2$ satisfying an appropriate differential equation and $g_{S^2}$ is the Fubini-Study metric on $S^2$ with area $\pi$.  In the patch $H_+$, $g_{S^2}$ is given by:
$$ g_{S^2} = \frac{|d\lambda_+|^2}{(1+|\lambda_+|^2)^2}
$$
The resolution parameter $a$  measures the size of the bolt $\C P^1$ ($a=0$ for the conifold).

\bigskip
The differential equation satisfied by the function $F_{a}$ is given
by imposing the Ricci-flat condition. It is convenient to write the equation in terms of $\gamma=r^2F'$, where the derivative is with respect to $r^2$. One gets:
$$\gamma'\gamma(\gamma+4a^2)=\frac{2}{3}r^2$$
where $\gamma'=\frac{d\gamma}{dr^2}$.
Integrating yields:
$$\gamma^3+6a^2\gamma^2-r^4=0$$ which gives the solution as in \cite{co}:
$$\gamma=-2a^2+4a^4N^{-\frac{1}{3}}+N^{\frac{1}{3}}, \
N(r)=\frac{1}{2}(r^4-16a^6+\sqrt{r^8-32a^6r^4})$$

\ni
The function $F_{a}'(r^2)$ is therefore given by:
\begin{equation} \label{F}
F_{a}'(r^2) = r^{-2} \l( -2a^2+4a^4N^{-\frac{1}{3}}+N^{\frac{1}{3}},\r) \,\, \mbox{where} \,
N(r)=\frac{1}{2}(r^4-16a^6+\sqrt{r^8-32a^6r^4})
\end{equation}

\ni For the conifold $Q_0$ with $a=0$, we have $F_{0}' =
r^{-\frac{2}{3}}$.

\bigskip

The K{\"a}hler form $\omega_{rc}(v,w)=g_{rc}(Jv,w)$ on the
resolved conifold can be expressed as a sum of two terms $
\omega_{rc} = d\,\alpha_{rc} + 4a^2\,\omega_{S^2} $, where the
one-form $\alpha_{rc}$ (compare \cite{im1}) is given by:

\begin{equation}\label{1-form}
\alpha_{rc}= F_{a}'(r^2)\ \Im \mbox{tr}(W^*dW) = F_{a}'(r^2)\Im
\l({\bar X}dX +{\bar Y}dY + {\bar U}dU + {\bar V}dV\r)
\end{equation}

\ni and $\omega_{S^2}$ is the standard K{\"a}hler form of area
$\pi$ on $S^2$. It can be expressed as
$\omega_{S^2}=d\alpha_{\pm}$ in the two patches $H_{\pm}$ on $S^2$
where:
$$\alpha_{\pm }=\frac{1}{2}\ \mbox{Im} \frac{\lambda_{\pm} d \bar
{\lambda_{\pm}}}{1+|\lambda_{\pm}|^2} \quad \mbox{on  } H_{\pm}$$
One can see that $\alpha_{rc}$ is invariant under the action of
$SO(4)= SO(4,\R)$. To show that $\alpha_{\pm}$ is also $SO(4)$-invariant, we
now compute the action of $SO(4)$ on $\C P^1$. The following three
matrices acting on the $z$-coordinates generate the $\sof(4)$ as a
Lie algebra:

\begin{equation}
B_1=\l(\begin{smallmatrix}0&1&0&0\\-1&0&0&0\\0&0&0&0\\0&0&0&0\end{smallmatrix}\r),
\
B_2=\l(\begin{smallmatrix}0&0&0&0\\0&0&0&0\\0&0&0&-1\\0&0&1&0\end{smallmatrix}\r),
\
A_3=\l(\begin{smallmatrix}0&0&0&0\\0&0&-1&0\\0&1&0&0\\0&0&0&0\end{smallmatrix}\r)
\end{equation}
Conjugating by the change of coordinate matrix $P$ (\ref{coord}),
we obtain the action on the $X,Y,U,V$-coordinates:

\begin{equation}
\tilde{B_1}=\l(\begin{smallmatrix}-i&0&0&0\\0&i&0&0\\0&0&0&0\\0&0&0&0\end{smallmatrix}\r),
\
\tilde{B_2}=\l(\begin{smallmatrix}0&0&0&0\\0&0&0&0\\0&0&i&0\\0&0&0&-i\end{smallmatrix}\r),
\
\tilde{A_3}=\frac{1}{2}\l(\begin{smallmatrix}0&0&-1&-1\\0&0&1&1\\1&-1&0&0\\1&-1&0&0\end{smallmatrix}\r)
\end{equation}

\ni In $H_+$, $\lambda_+=-\frac{X}{U}=-\frac{V}{Y}$ and so
$d\lambda_+=-\frac{d X}{U}+\frac{X dU}{U^2}$. The action of
$\tilde{B_1}$ on $\lambda_+$ is given by evaluating $d\lambda_+$
on the vector field corresponding to $\tilde{B_1}$. We obtain:
$\tilde{B_1}.\lambda_+=-i\lambda_+$. Similarly, on $H_ -$ we
compute that $\tilde{B_1}.\lambda_-=i\lambda_-$. The action of
$\tilde{B_2}$ is the same as for $\tilde{B_1}$ and the action of
$\tilde{A_3}$ is given by:
$$\tilde{A_3}.\lambda_+=\frac{1}{2}(1+\lambda_+^2), \
\tilde{A_3}.\lambda_-=-\frac{1}{2}(1+\lambda_-^2)$$ A computation
now shows that the form $\alpha_{\pm}$ is invariant under the
action of these three generators and hence under the
$SO(4)$-action. Therefore, the K{\"a}hler form and the metric are
$SO(4)$-invariant.

\ni The holomorphic volume $3$-form on $M$ has the form:

\begin{equation}\label{holo}
\Omega_{rc }= dU \wedge dY \wedge d\lambda_+ = dV \wedge dX \wedge d\lambda_{-}
\end{equation}

\ni in local coordinates and it is also invariant under the action
of $SO(4)$.

\bigskip
 We now want to describe cohomogeneity one special Lagrangian in the resolved conifold $M$. As in our previous
  paper \cite{im1}, the
 only subgroups of $SO(4)$ that we need to consider are the maximal torus $T^2$ and the subgroup $SO(3)$.
  We will use the same moment map techniques as in the previous paper. We recall the following
proposition which holds in general for any K{\"a}hler manifold
$M$.

\begin{proposition}\label{lagr}
 Let $G\subset SO(4)$ be a connected Lie subgroup with Lie
algebra $\g$ and moment map $\mu:M\to \g^*$ where $(M,
\omega_{rc})$ is the resolved conifold and $\mathcal O$ an orbit
of $G$ in $M$. Then the orbit is isotropic, i.e.
$\omega_{rc}|_{\mathcal O}\equiv 0$ if and only if ${\mathcal
O}\subseteq \mu^{-1}(c)$ for some $c\in Z(\g^*)$, where $Z(\g^*)$ is the centre of $\g$.
\end{proposition}

\bigskip

\section{$T^2$-invariant special Lagrangians}

Under the coordinate transform (\ref{action}), the action of the
maximal torus $T^2$ of $SO(4)$ on the $z$-coordinates, described
by the matrices:
$$\l\{\begin{pmatrix} \cos\theta_1 & -\sin\theta_1 & 0 &0 \\ \sin\theta_1 &\cos\theta_1 & 0 & 0\\ 0& 0 &
\cos\theta_2 & -\sin\theta_2 \\ 0& 0& \sin\theta_2 & \cos\theta_2
\\\end{pmatrix} \r\}$$
correspond to the following matrices with respect to the
$(X,Y,U,V)$-coordinates:

$$\l\{\begin{pmatrix} e^{-i\theta_1} & 0 & 0 &0 \\ 0 & e^{i\theta_1} & 0 & 0\\ 0& 0 & e^{i\theta_2}
& 0 \\ 0& 0& 0 & e^{-i\theta_2} \end{pmatrix}\r\}$$

\bigskip
\ni The $T^2$-action on the patch $H_+$ with coordinates
$(U,Y,\lambda)$ is given by: $$g.(U,Y,\lambda_+)=(e^{i\theta_2}U,
e^{i\theta_1}Y,e^{-i(\theta_1+\theta_2)}\lambda_+)$$ and on the
patch $H_-$ with coordinates $(X,V,\mu)$ is given by:
$$g.(X,V,\lambda_{-})=(e^{-i\theta_1}X,
e^{-i\theta_2}V,e^{i(\theta_1+\theta_2)}\lambda_{-})$$ The
following result gives all the special Lagrangian $3$-folds of the
resolved conifold invariant under the action of $T^2$.

\begin{theorem}
The special Lagrangian submanifolds in the resolved conifold $M$
with coordinates $(X,Y,U,V,\lambda_1,\lambda_2)$ described by
equations (\ref{rescon}) with the Calabi-Yau structure, which are
invariant under the action of the maximal torus $T^2$ of $SO(4)$
are given by the equations:

\begin{align}
&\frac{1}{2} F_{a}'(r^2)(|X|^2-|Y|^2)+4a^2\frac{|\lambda_2|^2}{|\lambda_1|^2+ |\lambda_2|^2} = c_1 \notag\\
& \frac{1}{2}F_{a}'(r^2)(|V|^2-|U|^2)+4a^2 \frac{|\lambda_2|^2}{|\lambda_1|^2+ |\lambda_2|^2}= c_2 \label{eqt2}\\
&\Im (XY)= c_3\notag
\end{align}

\ni where $F_{a}'$ is given by (\ref{F}) and $c_1,c_2$ and $c_3$
are real constants. \label{t}
\end{theorem}

\bigskip
\ni
\pf
\ni
Using the two infinitesimal generators:

$${\tilde B_1} =\l(\begin{smallmatrix}-i&0&0&0 \\
0& i&0&0\\0&0&0&0\\0&0&0&0\end{smallmatrix}\r), \
{\tilde B_2} =\l(\begin{smallmatrix}0&0&0&0 \\
0& 0&0&0\\0&0&i&0\\0&0&0&-i\end{smallmatrix}\r)$$

\ni
of the $T^2$-action on the resolved conifold, the moment map is given by:

\begin{align*}
&\mu: M \to ({\mathfrak t}^2)^*\simeq \R^2, \\
 &\mu =
\l(\frac{1}{2} F_{a}'(r^2)(|X|^2 - |Y|^2) +
4a^2\,\mu_{S^2}\,, \,\,\frac{1}{2} F_{a}'(r^2)(|V|^2 -
|U|^2) + 4a^2\,\mu_{S^2}\r)
\end{align*}
where
$\mu_{S^2}(\lambda_1,\lambda_2)=\frac{|\lambda_2|^2}{|\lambda_1|^2+|\lambda_2|^2}$
is the moment map of the standard $S^1$-action on $\C P^1$.

\medskip
Since $({\mathfrak t}^2)^*=\R^2$ is abelian, it follows from
Proposition $\ref{lagr}$ that any $T^2$-invariant special
Lagrangian $3$-fold $L$ in $Q^3$ lies in a level set
$\mu^{-1}(c)$, where $c=(c_1,c_2)\in \R^2$. The first two
equations enforce this condition and ensures that the submanifold
is Lagrangian. In order
to impose the special Lagrangian condition at a given point $p=
(U,Y,\lambda_+)$, we compute $\Omega_{rc}$ on the coordinate patch
$H_+$ with coordinates on the three tangent vectors $Y_1={\tilde
B_1}\,p = (0,iY,-i\lambda_+), \ Y_2={\tilde B_2}\,p =
(iU,0,-i\lambda_+)$ and $Y_3= \dot{p}=
(\dot{U},\dot{Y},\dot{\lambda}_+)$:

\begin{align}
\Omega_{rc}(Y_1,Y_2,Y_3) &=(dU \wedge dY \wedge d
\lambda_+)(Y_1,Y_2,Y_3) \notag\\
 &=  \l| \begin{matrix}
0           &  iU             & \dot{U} \\
iY          &   0             & \dot{Y} \\
-i\lambda_+ &  - i\lambda_+    & \dot{\lambda}_+
\end{matrix} \r|\\
&= \dot{U}\,Y\,\lambda_+ + U\,\dot{Y}\,\lambda_+ +
U\,Y\,\dot{\lambda}_+=(UY\lambda_+)^\cdot\notag
\end{align}

\bigskip
\ni Integrating the condition $\Im \, \Omega_{rc} = 0$, we obtain
$ \Im(UY\lambda_+)= c$ for $c \in \R$. Using $\lambda_+ = -
\frac{X}{U}$ we finally obtain the third equation of the theorem.
Note that in the other coordinate patch $H_ -$ we will obtain $\Im
(VX\lambda_ -)=$ c which is the same equation since $\lambda_
-=-\frac{Y}{V}$.

\hfill $\Box$

\bigskip

\ni {\it Remark 1:} Equations (\ref{eqt2}) are obviously
$T^2$-invariant and linearly independent at a generic point The above
family of $T^2$-invariant special Lagrangian $3$-folds foliates the
resolved conifold $ M =  {\mathcal O}(-1)\oplus {\mathcal O}(-1)$. 
The generic orbit is $T^2 \times \R$
where $T^2$ is an orbit of the maximal torus in $SO(4)$. For the
special values $c_1=c_2, \ c_3=0$, these special Lagrangians
intersect the $\C P^1$ in circles.

\smallskip

\ni {\it Remark 2:} Asymptotically, as $r^2 \to \infty$ or
equivalently as $a \to 0$, these special Lagrangians approach the
special Lagrangian cone on $T^2$ in the conifold described by the
equations:
\begin{align}
&|X|^2-|Y|^2 = 0 \notag\\
&|V|^2-|U|^2 = 0\notag  \\
&\Im (XY)= 0 \notag
\end{align}
which is the same asymptotic cone that was described in
\cite{im1}. In $z$ coordinates on the conifold, the equations
become: $\Im(z_0\bar{z_1})=\Im(z_2\bar{z_3})= \Im (z_0^2+z_1^2)=
0$


\section{$SO(3)$-invariant special Lagrangians}

Let $G=SO(3)$ be the subgroup of $SO(4)$ embedded as:
$$\l\{\begin{pmatrix}1&0 \\ 0& A\end{pmatrix},\ A\in SO(3)\r\}$$
The $SO(3)$-invariant special Lagrangian $3$-folds in the resolved
conifold $M$ are characterized in the following theorem:

\begin{theorem}\label{thm1}
The $SO(3)$-invariant special Lagrangian submanifolds in the
resolved conifold $M$, with coordinates
$(X,Y,U,V,\lambda_1,\lambda_2)$ satisfying equation (\ref{coneq})
are given by:
$$L_c=\{\tilde{g}.p\in M \,| \,  g\in SO(3),\ p=(0,Y,0,0), Y \neq 0,  \ {\mbox Re } (Y^2)=c\}$$
where $\tilde{g}$ is the conjugated action on the $(X,Y,U,V)$-
coordinates.
\end{theorem}
\pf

\ni Let $ A_1=\l(\begin{smallmatrix}0&0&0&0 \\
0& 0&0&0\\0&0&0&-1\\0&0&1&0\end{smallmatrix}\r), \,
A_2=\l(\begin{smallmatrix}0&0&0&0 \\
0& 0&0&1\\0&0&0&0\\0&-1&0&0\end{smallmatrix}\r), \,
A_3=\l(\begin{smallmatrix}0&0&0&0 \\
0& 0&-1&0\\0&1&0&0\\0&0&0&0\end{smallmatrix}\r)$ be the
infinitesimal generators of $SO(3)$. Using the change of
coordinate matrix $P$ (\ref{coord}), the action on the
$(X,Y,U,V)$-variables is given by:

$${\tilde A_1}=\l(\begin{smallmatrix}
0  &  0  &  0 &  0   \\
0  &  0  &  0  &  0   \\
0  &  0  &  i&  0  \\
0 &0   &  0 &  -i
\end{smallmatrix}\r),\
{\tilde A_2}= \frac{1}{2}\l(\begin{smallmatrix}
0  &  0  &  -i &  i   \\
0  &  0  &  i  &  -i   \\
-i  &  i   &  0 &  0  \\
i   & -i   &  0 &  0
\end{smallmatrix}\r),\
{\tilde A_3} = \frac{1}{2}\l(\begin{smallmatrix}
0   &  0   &   -1  &  -1  \\
0   &  0   &    1   &  1   \\
1   &  -1  &    0  &   0   \\
1   &   -1 &   0   &   0
\end{smallmatrix}\r)$$

\ni The moment map for these generators is $\mu:M\to \sof (3)^*$ given by:\\

$$\mu = \l(\frac{1}{2}(|U|^2-|V|^2),\
\frac{1}{2}\Im((U-V)(\bar{Y}-\bar{X})),
\frac{1}{2}\Im((Y-X)(\bar{U}+\bar{V}))\r)+\mu_{S^2}$$ where
$$\mu_{S^2}(\lambda_+)=\l(
\frac{|\lambda_+|^2}{1+{|\lambda_+|^2}},-\frac{1}{2(1+|\lambda_+|^2)}\Re(\lambda_+
 -|\lambda_+|^2\bar{\lambda_+}),-\frac{1}{2(1+|\lambda_+|^2)}\Im(\lambda_+
 +|\lambda_+|^2\bar{\lambda_+})  \r)$$
on the patch coordinate patch $H_+$. A similar expression can be
computed on $H_ -$.

Let $p= (X, Y, U, V) \in M$. Since $Z(\sof(3)^*)={0}$, we need to
start at a point in $\mu^{-1}(0)$. We make the simplifying assumption that
$U=V=0$. Using the fact that $XY=UV$, we
can assume that $X=0$ or $Y=0$. If we choose $Y=0$, then the point
$(X,0,0,0)\not\in \mu^{-1}(0)$ since $\mu(X,0,0,0)$ has the first
component equal to 1. Hence we must choose $p= (0, Y, 0, 0)$, where $Y \neq 0$.
Since $\lambda_+=-\frac{V}{Y}=0$, $p\in\mu^{-1}(0)$. We now look
for curves $Y(s)$ in the complex $Y$-plane, which after applying
the $SO(3)$-action give rise to special Lagrangians of the form
$L=\tilde{g}.Y(s)$, $g\in SO(3)$.

\ni
The tangent plane at $p$ to $L$ is spanned by the vectors:
$\{v_1= {\tilde A_2}\,p \, , \, v_2= {\tilde A_3}\,p \, , \, v_3= \dot{p} \, \}$, 
where the dot denotes differentiation with respect to the parameter $s$.
The generic orbit of the $SO(3)$-action is an $S^2$, the three infinitesimals generators are
linearly dependent (Note that ${\tilde A_1}\,p=0$). The tangent space of
$L$ at the point $p$ is therefore spanned by the vectors:
$$\l\{v_1=\frac{iY}{2}\begin{pmatrix} 0\\0\\1\\-1\end{pmatrix},\
v_2=\frac{Y}{2}\begin{pmatrix} 0\\0\\
-1\\ -1 \end{pmatrix} ,\ v_3=\begin{pmatrix}
0\\\dot{Y}\\0\\0\end{pmatrix}\r\}$$ \ni $L$ is invariant under the
$SO(3)$-action and $\omega_{rc}\mid_L=0$, since $L$ lies in the
zero level set of the moment map. We now impose the special
Lagrangian condition $\mbox{ Im } \Omega_{rc}\mid_L=0$. Working
on the patch $H_ +$ with coordinates $(U,Y,\lambda_+)$, the
holomorphic volume form is given by equation (\ref{holo}):
$$\Omega_{rc}=dU\wedge dY\wedge d\lambda_{+}$$
Since $\lambda+-=-\frac{V}{Y}$, we can instead use coordinates
$(U,Y,V)$ on  and the holomorphic volume form in these coordinates
becomes:
$$\Omega_{rc}=dU\wedge dY\wedge d\lambda_{+}=-dU\wedge dY\wedge \frac{VdY-YdV}{Y^2}=-\frac{1}{Y}dU\wedge dY\wedge dV$$

\ni Therefore, $$\Omega_{rc}(v_1,v_2,v_3)=-\frac{1}{Y}dU\wedge
dY\wedge dV(v_1,v_2,v_3)$$ Computing yields:
$$ -\frac{1}{Y}\l| \begin{matrix}
\frac{iY}{2}&-\frac{Y}{2}&0\\
0 & 0 & \dot{Y} \\
-\frac{iY}{2} & -\frac{Y}{2}&0
\end{matrix} \r|
=  -\frac{i}{2}Y\dot{Y}
$$
\ni Integrating, the condition $\mbox {Im } \Omega_{rc}=0$
becomes:
\begin{equation}\label{so3} \mbox{Re} (Y^2)=c
\end{equation}
 \ni
where $c$ is a real constant.
\ni Let $Y=u+iv$. Then equation (\ref{so3}), gives a family of
hyperbolas in the $(u,v)$-plane:
$$u^2-v^2=c$$

\ni
The special Lagrangian $L_c$ is topologically $S^2\times \R$
and has two components, each asymptotic to the special Lagrangian
cones on $S^2$ given by setting $c=0$ (obtained from the lines
$u+v=0$ and $u-v=0$) in the conifold. This case is reminiscent of
the flat case studied in \cite{hl}, but of a different dimension.

\medskip
\ni 
{\it Remark 1:} 
$L_c$ does not intersect the zero-section (the
bolt) $S^2$ of the resolved conifold, since $r^2 = |Y|^2$ is preserved under the $SO(3)$.

\smallskip
\ni 
{\it Remark 2:} 
 Note that asymptotically,
the $SO(3)$-invariant special Lagrangian in the resolved conifold
approach the same special Lagrangian cones in the conifold as the
$SO(3)$-invariant special Lagrangian in the deformed conifold that
we found in our previous paper \cite{im1}.

\pagebreak

\bigskip
\bigskip
\bigskip

\ni 
{\it Department of Mathematics \& Statistics, McMaster University \\
e-mail address: ionelm@@math.mcmaster.ca} \\

\medskip
\ni 
{\it Department of Mathematics \& Statistics, McMaster University\\
e-mail address: minoo@@mcmaster.ca \\}

\end{document}